\documentclass[11pt]{amsart}

\usepackage{amsmath, amssymb, amsfonts}
\usepackage{graphicx, overpic}
\usepackage{fullpage, multirow, tabulary}
\usepackage{color}
\usepackage{enumitem}  
\usepackage[normalem]{ulem}
\usepackage[colorlinks=true, pdfstartview=FitV, linkcolor=blue, citecolor=blue, urlcolor=blue]{hyperref}
\usepackage{array, scalerel}

\definecolor{cerulean}{rgb}{0,.48,.65} 
\definecolor{magenta}{rgb}{.5,0,.5} 
\definecolor{dred}{rgb}{.5,0,0} 
\definecolor{green}{rgb}{0,.5,0} 
\definecolor{blue}{rgb}{0,0,0.5} 
\definecolor{black}{rgb}{0,0,0} 
\definecolor{dgreen}{rgb}{0,.3,0} 
\definecolor{vdred}{rgb}{.3,0,0} 
\definecolor{red}{rgb}{1,0,0} 
\definecolor{salmon}{rgb}{0.98,0.50,0.45} 
\definecolor{gray}{rgb}{.5,.5,.5} 
\definecolor{seagreen}{rgb}{0.13,0.70,0.67} 
\definecolor{chartreuse}{rgb}{0.40,0.80,0.00}
\definecolor{cornflower}{rgb}{0.39,0.58,0.93} 
\definecolor{gold}{rgb}{0.60,0.48,0.00}
\definecolor{orange}{rgb}{1,0.25,0.05} 
\definecolor{bcyan}{rgb}{0.15,0.59,0.75}

\theoremstyle{plain}
\newtheorem{thm}{Theorem}[section]
\newtheorem{lemma}[thm]{Lemma}
\newtheorem{prop}[thm]{Proposition}

\newtheorem{corollary}[thm]{Corollary}

\numberwithin{equation}{section}
\numberwithin{figure}{section}
\numberwithin{table}{section}

\theoremstyle{definition}

\newtheorem{remark}[thm]{Remark}

\newcommand{\R}{\mathcal{R}}
\renewcommand{\S}{\mathcal{S}}
\newcommand{\U}{\mathcal{U}}
\newcommand{\V}{\mathcal{V}}
\newcommand{\N}{\mathbb{N}}

\newcommand{\X}{\mathcal{X}}
\newcommand{\wo}{\overline{w}}

\newcommand{\set}[1]{\left\{#1\right\}}

\def\N{\mathbb{N}}

\def\P{\mathcal{P}}
\def\Q{\mathcal{Q}}

\def\Dist{\hbox{\rm Dist}}

\def\Diam{\hbox{\rm Diam}}

\def\Area{\hbox{\rm Area}}
\def\Rad{\hbox{\rm Rad}}

\def\onto{{\kern3pt\to\kern-8pt\to\kern3pt}}

\newcommand{\hide}[1]{}

\title{The distortion of finitely generated subgroups \\ in hyperbolic groups}

\author{Pallavi Dani and Timothy Riley}
\date{\today}

\thanks{We gratefully acknowledge the financial support of the National Science Foundation, NSF DMS-2407104  (PD), NSF GCR-2428489 (TR)} 

\thanks{\emph{ORCID:} 0000-0001-7653-2618 (PD),  0009-0004-3699-0322 (TR)}   



\begin{document}

\begin{abstract}
We show that the range of distortion functions of finitely generated subgroups of hyperbolic groups is essentially as broad as the range of time-complexity functions $T(n)$ of Turing machines. More precisely, we show that for each such $T(n) \geq n$ for which  $T(n)^3$ is superadditive, there is a hyperbolic group $G$ with a finitely generated normal subgroup $N$ such that $\Dist_N^G(n) \simeq 2^{T(n)^3}$. Our approach is to estimate the distortion of the kernels in Wise's HNN version of the Rips construction, and then to look at what these estimates yield when the $S$-machine groups of Sapir, Birget, and Rips are the inputs to the construction. 
\end{abstract}

\maketitle

\section{Introduction}\label{sec:intro}

\subsection{Our results}

In the setting of finitely generated groups, the distortion function quantifies how a subgroup sits metrically inside a group. This article explores which distortion functions can occur for subgroups of hyperbolic groups by relating the range of possibilities to the time-complexity functions of Turing machines.

\begin{thm} \label{thm: time to distortion}
Let $T(n) \geq n$ be the time function of a (not necessarily deterministic) Turing machine.  Assume that $T(n)^3$ is superadditive.  Then there exist a hyperbolic group $G$ and a finitely generated normal subgroup $N$ of $G$ such that the distortion function for $N$ in $G$ satisfies
$$\Dist_N^G(n) \simeq 2^{T(n)^3}.$$
\end{thm}

In Section~\ref{sec:preliminaries}, we provide some preliminaries, including the definitions of $\simeq$ and $\Dist_N^G$ used here.

As we will explain in Section~\ref{sec:deductions}, the above theorem is related to the study of Dehn functions in \cite{SBR}. The same is true of our next theorem, which identifies a more explicit sweep of distortion functions for subgroups of hyperbolic groups:

\begin{thm} \label{thm: examples}
Let $\beta \ge 3$ be a real number for which there exists $C>0$ such that for every $m$, the first $m$ digits of $\beta$ can be computed in time $\leq C2^{2^{Cm}}$. Then there exist a hyperbolic group $G$ and a finitely generated normal subgroup $N$ of $G$ such that 
$$ \Dist_N^G(n) \simeq 2^{n^{\beta}}.$$
\end{thm} 

Examples of numbers $\beta \geq 3$ satisfying the hypothesis here include rational numbers, $\pi$, $e+1$, and $e\pi$ \cite{Brent} (see also \cite[Section~1]{SBR}).

Theorems~\ref{thm: time to distortion} and \ref{thm: examples} show that all \emph{reasonable} functions growing at least like $n \mapsto 2^{n^3}$ are distortion functions of finitely generated normal subgroups of hyperbolic groups. We recently proved that, for every pair of integers $p > q > 0$, the function $2^{n^{p/q}}$ is the distortion function of a finite-rank free subgroup of a hyperbolic group \cite{DR}. Our results here give a more detailed account of which distortion functions can occur, but above $2^{n^3}$ rather than above $2^n$, and at the expense of losing finite presentability. 

That our theorems yield no   subgroups of subexponential distortion is to be expected. Kapovich, following Gromov, proved that if the distortion of a finitely generated subgroup of a hyperbolic group is subexponential, then that subgroup is undistorted, meaning that its distortion function is $\simeq n$ \cite[Proposition~2.6]{Ksubexp}.

A more detailed history of subgroup distortion in hyperbolic groups can be found in \cite{DR}.

\subsection{Our strategy}

The Rips construction, $1 \to N \to G \to Q \to 1$, associates a hyperbolic group $G$ to any finitely presented group $Q$ \cite{Rips}.  We will prove inequalities relating the distortion of the Rips kernel $N$ in $G$ to the Dehn and isodiametric functions of $Q$.
We will then take $Q$ from the family of finitely presented groups constructed by Sapir, Birget, and Rips~\cite{SBR}, which realize a wide range of Dehn functions, and find that the scope of the resulting distortion functions is comparable to that of the time-complexity functions of Turing machines.

However, instead of Rips's original construction, we will use Wise's HNN version \cite{Wise1}. The additional structure it affords $G$ allows us to prove (in Theorem~\ref{Main}) upper and lower bounds on the distortion of $N$ in $G$ that are sufficiently tight for us to be able to deduce Theorems~\ref{thm: time to distortion} and~\ref{thm: examples}.

These bounds will come from an analysis of the geometry of van~Kampen diagrams over $G$. In Wise's construction, these diagrams consist of corridors, which we will examine, for example in Lemmas~\ref{t corridors} and \ref{key corridor lemma}, using techniques from our recent paper \cite{DR}. Small-cancellation arguments will imply that, in reduced van~Kampen diagrams, no pair of these corridors can stay close for more than a constant distance. It will follow that the corridors nest in a doubling pattern with the depth of this nesting  essentially being the diameter of the image of the diagram in $Q$. This doubling leads (in Theorem~\ref{Main}) to a relationship between the distortion of $N$ in $G$, the area of diagrams in $Q$, and an exponential of the  diameters of these diagrams.

\subsection{Wise's HNN version of the Rips construction} \label{sec: Wise HNN}
Given a finitely presented group
$$Q = \langle a_1, \dots, a_m \mid r_1, \dots, r_k \rangle,$$
Wise defines a group $G$ by setting $\S = \set{a_1, \dots, a_m, x_1, x_2, t}$ and taking
\begin{equation} \label{Rips presentation}
G   =   \scaleleftright[1.9ex]{\bigg\langle}{    \   \S \   \rule{0mm}{12mm}  \left|  \ \rule{0mm}{12mm} \parbox{97mm}{$r_i = X_i \, t \, X'_i \,t^{-1} $ for $1\le i\le k,$   \\[6pt] \ $a_j^\epsilon x_l a_j^{-\epsilon} = X_{\epsilon jl} \,t \,X'_{\epsilon jl} \, t^{-1} $ \  for \ $1\le j\le m$, \  $\epsilon =\pm 1$, \ $l= 1, 2$ \\[6pt] $a_j^\epsilon t a_j^{-\epsilon} = X_{\epsilon j}\, t \, X'_{\epsilon j}$  \  for \ $1\le j \le m$, \  $\epsilon =\pm 1$  } \right. }{\bigg\rangle}
\end{equation}
where the $2k+12m$ variously decorated symbols $X$ appearing in the defining relations form a set of words $\X$ that are pairwise disjoint subwords of the infinite Rips word $x_1 x_2^1 \, x_1 x_2^2 \, x_1 x_2^3 \cdots$. We will specify these subwords in Section~\ref{sec: small cancellation}; they are chosen so that the set of relations satisfies small-cancellation conditions implying that $G$ is hyperbolic and, indeed, $\textup{CAT}(-1)$, cubical, virtually special, and residually finite.

Let
\begin{equation}\label{eq:UV}
\U = \set{U_i \mid 1 \leq i \leq k+6m} \qquad \text{and} \qquad \V = \set{V_i \mid 1 \leq i \leq k+6m}
\end{equation}
be sets of words such that the $k+6m$ defining relations in \eqref{Rips presentation} can be expressed as $t U_i t^{-1} = V_i$.
In Section~\ref{sec: small cancellation}, we will explain how small-cancellation conditions imply that $\U$ and $\V$ are free bases for subgroups of the free group $F = F(a_1, \ldots, a_m, x_1, x_2)$. It will follow that $G$ is an HNN extension of $F$ with stable letter $t$ (see Corollary~\ref{cor:hnn structure}).

Let $\R$ be the set of the $k+6m$ words $\sigma \tau^{-1}$ on $\S^{\pm 1}$ obtained from the defining relations $\sigma = \tau$ in \eqref{Rips presentation}. Thus $G = \langle \S \mid \R \rangle$.

Killing $t$, $x_1$, and $x_2$ defines an epimorphism from $G$ onto $Q$ whose kernel is the subgroup $N$ generated by $t$, $x_1$, and $x_2$. By \cite{Bieri}, $N$ is finitely presentable if and only if $Q$ is finite. %

\subsection{Relating area and diameter to distortion}

An observation concerning the Rips construction, made in \cite[Section~3.4]{AO}, \cite[Corollary~8.2]{Farb}, \cite[Section~3, $3.K_3''$]{Gromov}, and \cite{PittetThesis},  is that if $Q$ is a finitely presented group with unsolvable word problem, then the membership problem for $N$ in $G$ is undecidable.

These undecidability statements are reflected in functions that quantify the corresponding decision problems. If $Q$ has unsolvable word problem, then its Dehn function cannot be bounded above by a recursive function: otherwise, an exhaustive process of applying relators would solve the word problem. Similarly, if the membership problem for $N$ in $G$ is undecidable, then the distortion function $\Dist^G_N$ cannot be bounded above by a recursive function: otherwise, since $G$ is hyperbolic and therefore has decidable word problem, an exhaustive search would solve the membership problem for $N$ in $G$.

This suggests that there should be a general relationship between the Dehn function of $Q$ and the distortion function of $N$ in $G$, and indeed we will give one in Theorem~\ref{Main} below. However, the relationship we will present also involves the diameter (isodiametric) function.  The diameter function is always dominated by the Dehn function, and the Dehn function is bounded above by a double exponential of its diameter function \cite{Cohen}. So the diameter function of a finitely presented group is bounded from above by a recursive function (indeed, is a recursive function) if and only if the same is true of the Dehn (area) function, and a bound on distortion in terms of one of these functions  would imply a bound in terms of the other.  However, combining Dehn and diameter functions appropriately will allow us to  give tighter upper and lower bounds on $\Dist_N^G (n)$ than we could with one of these functions alone.
To that end, we will require area and diameter bounds that can be realized simultaneously: an \emph{area-diameter pair} $(f,g)$ for a finitely presented group $Q$ is a pair of non-decreasing functions $f,g \colon \N \to \N$ such that every word $w$ of length at most $n$ representing the identity in $Q$ admits a van~Kampen diagram with area at most $f(n)$ and diameter at most $g(n)$. We define this in detail in Section~\ref{sec:preliminaries}.

\begin{thm} \label{Main}
Let $1 \to N \to G \to Q \to 1$ be the exact sequence obtained by applying the Rips--Wise construction to an infinite finitely presented group $Q$, as described in Section~\ref{sec: Wise HNN}. Suppose $(f,g)$ is an area-diameter pair for $Q$. 
Then 
there exist constants $K_1, K_2 >1$ such that 
the distortion of $N$ in $G$ satisfies
\begin{equation}\label{distortion bounds}
K_1^{\left(\Diam_Q(n)-\lfloor \frac n 2 \rfloor\right)} \preceq \Dist_N^G (n) \preceq (f(n) +n) \  K_2^{g(n)}.
\end{equation}
\end{thm}

We will prove the upper bound on $\Dist_N^G (n)$ in Section~\ref{sec: upper bound} and the lower bound in Section~\ref{sec: lower bound}.

\subsection{Deducing Theorems~\ref{thm: time to distortion} and \ref{thm: examples} from Theorem~\ref{Main}} \label{sec:deductions}

Sapir, Birget, and Rips proved in \cite[Theorem~1.3]{SBR} that if $\mathcal{M}$ is a (not necessarily deterministic) Turing machine with time function $T(n)$ such that $T(n)^4$  is superadditive, then there exists a finitely presented group $G(\mathcal{M})$ whose Dehn function is $\simeq T(n)^4$ and whose smallest isodiametric function is $\simeq T(n)^3$.  
Moreover, Proposition~11.1 of \cite{SBR}  makes clear that these two functions are realized by the same van~Kampen diagrams (see also the closing paragraphs of Section~1 and the lower bound estimates of Propositions~12.1 and 12.2).  So \cite{SBR} provides an area-diameter pair $(f,g)$ for $G(\mathcal{M})$ with $f(n) \simeq T(n)^4$ and $g(n) \simeq T(n)^3$. Taking $Q = G(\mathcal{M})$ in Theorem~\ref{Main} therefore yields Theorem~\ref{thm: time to distortion} --- the superadditivity of $T(n)^3$ implies superadditivity of $T(n)^4$ (because $T(n)$ is  non-decreasing by definition) and allows the constants implicit in the $\preceq$-relation to be pulled through --- and Lemma~\ref{lem: exp simeq}, applied twice, lets us replace the constants $K_1$ and $K_2$ by $2$.

As an application of their methods, Sapir, Birget, and Rips produce, for every real number $\alpha \geq 4$ whose decimal expansion is computable in the sense stated in Theorem~\ref{thm: examples}, a finitely presented group whose Dehn function grows as $\simeq n^{\alpha}$ and whose diameter function (the \emph{smallest isodiametric function} in their terminology) grows as $\simeq n^{3\alpha/4}$ \cite[Corollary~1.4]{SBR}.  As before, these bounds are simultaneously realizable. We may therefore apply Theorem~\ref{Main}. Set $\beta = 3\alpha/4$. The assumption on the computability of the decimal expansion of $\alpha$ follows from the same assumption on that of $\beta$, and, up to $\simeq$, the upper and lower bounds of Theorem~\ref{Main} both simplify to $2^{n^{\beta}}$ (in light of Lemma~\ref{lem: exp simeq}).  So we have Theorem~\ref{thm: examples}.

\subsection{Statement about AI usage} We conceived and wrote this paper without AI assistance. We then consulted Claude and Gemini for proofreading, resulting in a handful of minor corrections.

\section{Preliminaries}  \label{sec:preliminaries}

\subsection{Words}  
 
A \emph{word $w$ on} a set of letters $S$ is an expression $s_1^{\varepsilon_1} \cdots s_l^{\varepsilon_l}$     where $l \geq 0$, $s_i \in S$, and $\varepsilon_i = \pm 1$ for all $i$.  
Its length  $l$ is denoted $|w|$.   The word metric $d_S(g, h)$ on a group $G$ generated by $S$ gives the length of a shortest word on $S$ that represents $g^{-1}h$.

\subsection{Equivalent growth rates}  \label{sec: growth rates}

For functions   $f$ and $g$ mapping  $\mathbb{N} \to \mathbb{N}$ or $\mathbb{R}_{\geq 0} \to \mathbb{R}_{\geq 0}$ write $f \preceq g$ when there exists an integer $C>0$ such that $f(n) \leq Cg(  Cn+C ) +Cn+C$ for all $n \geq 0$ in the domain, and $f  \simeq g$ when $f \preceq g$ and $g \preceq f$.  

\begin{lemma} \label{lem: exp simeq}
Suppose $\phi$ and $\psi$ are non-decreasing functions $\N \to \N$ such that  $\phi$ grows at least linearly (meaning there exists $\zeta>0$ such that $\phi(n) \geq \zeta n$ for all $n$) and $\psi$ is  superadditive (meaning $\psi(x + y) \geq \psi(x) + \psi(y)$ for all $x,y \in \N$).  If $\phi \simeq \psi$, then $K^{\phi(n)} \simeq 2^{\phi(n)}$ for all $K >1$.    
\end{lemma}

\begin{proof} 
We claim that 
for all $A>0$, there exists an integer $M \ge 1$ such that  for all  $n \in \N$,    
\begin{equation} \label{phi bound}
A \, \phi(n) \leq \phi(M n + M) + M n + M.
\end{equation}
Indeed, choose an integer $C \geq 1$ such that
\begin{equation} \label{back and forth}
\phi(n)\leq C\psi(Cn+C)+Cn+C \ \ \text{ and }  \ \  
\psi(n)\leq C\phi(Cn+C)+Cn+C,
\end{equation}
put $x=Cn+C$ so that by superadditivity $L\psi(x) \leq \psi(Lx)$ for all integers $L \geq 1$, and then: 
$$
\begin{aligned}
A\phi(n)
& \leq A C \psi(x) + A(Cn + C)\\
& \leq \frac{A C}{L} \psi(Lx) + A(Cn + C)\\
& \leq \frac{AC^2}{L} \phi(C L x + C)
       +\frac{AC}{L}(Lx+C)+A(Cn+C).
\end{aligned}
$$
Take $L\geq AC^2$, so that the coefficient of $\phi(C L x + C)$ is then at most one. Since $x=Cn+C$, all the remaining terms are bounded above by a linear function of $n$.
Finally, since $\phi$ is nondecreasing, we can take $M$ sufficiently large that
$C L x + C \leq M n + M$ and that the linear error is at most $Mn + M$,  proving \eqref{phi bound}.

Per hypothesis, let $\zeta>0$ satisfy
$\phi(n)\geq\zeta n$ for all $n$.  
Then by \eqref{phi bound}, 
$$A\phi(n)
\leq  (1+\zeta^{-1})\phi(Mn+M).$$
So, for all $a>0$, we can choose $A=a(1+\zeta^{-1})$ and obtain
$a\phi(n)\leq\phi(Mn+M)$.
Taking $a=\log_2 K$, we get
$$K^{\phi(n)} =2^{(\log_2 K)\phi(n)} \leq 2^{\phi(Mn+M)}$$
and so $K^{\phi(n)}\preceq 2^{\phi(n)}$.

Interchanging the roles of $K$ and $2$ gives the reverse inequality. Hence
$K^{\phi(n)} \simeq 2^{\phi(n)}$.
\end{proof}  

\subsection{Distortion}  
  
For a group $G$ with finite generating set $S$ and a subgroup $H$ with finite generating set $T$, the \emph{distortion  function}  $\Dist^{G}_{H} : \N \to \N$ quantifies how $H$ sits metrically in $G$:
$$\Dist^{G}_{H}(n)   :=    \max \set{  \  d_T(1,g)  \  \mid  \   g \in {H} \textup{ with }  d_S(1,g)  \leq n   \ }.$$
 Other choices of finite generating sets for $G$ and $H$ yield   distortion functions that are  $\simeq$-equivalent.

\subsection{Dehn, diameter, and radius  functions} \label{sec: Dehn etc}

We refer to  \cite{Bridson6, RileyNotes, Sapir} for more detailed accounts of the following.

A \emph{presentation} $\langle S \mid R \rangle$ comprises an alphabet $S$ (\emph{generators}) and a set $R$   (\emph{defining relators}) of words on $S$, and defines the group $G = F(S) / \langle \! \langle R \rangle \! \rangle$, the quotient of the free group $F(S)$ by the smallest normal subgroup  containing the set of words $R$.   So a word $w$ on $S$ represents $1$ in $G$ if and only if 
\begin{equation} \label{prod}
w    =    \left( {u_1}^{-1} r_1^{\epsilon_1} u_1 \right) \cdots  \left( {u_{\ell}}^{-1} r_{\ell}^{\epsilon_{\ell}} u_{\ell} \right) \ \text{ in }  F(S)  
\end{equation}  for some ${\ell} \in \N$, some $\epsilon_i = \pm 1$, some words $u_i$ on $S $, and some $r_i \in R$.   

A presentation expressed as $\Big\langle S  \,  \left| \rule{0mm}{4mm} \right. \,  \set{v_i = w_i}_{i \in I} \Big\rangle$ is shorthand for $\Big \langle S  \,  \left| \rule{0mm}{4mm} \right. \,  \set{v_i   w_i^{-1}}_{i \in I} \Big\rangle$.   

A \emph{van~Kampen diagram for $w$ over $G$}  is a finite planar contractible 2-complex $\Delta$ such that each edge is oriented and is labelled by a generator so that around each face one reads a defining relator, and around the diagram's boundary circuit one reads $w$ starting from a base vertex $\ast$. 

Such a $\Delta$ demonstrates how $w=1$ in $G$ follows from the defining relators. Indeed, \emph{van~Kampen's Lemma} \cite{vK} states that a word $w$ represents the identity in a group $G$ with finite presentation $P$ if and only if it admits a van~Kampen diagram over $P$. In outline, a \emph{lollipop diagram} $\Delta'$ whose boundary word is the right-hand side of \eqref{prod} can be folded to form a van~Kampen diagram for $w$. Each $u_i$ in $\Delta'$ labels an edge-path $\rho_i$ from a  face whose boundary circuit is labelled by $r_i$ (with one orientation or the other) to a basepoint $\ast$. Conversely, a van~Kampen diagram for $w$ can be \emph{cut open} to form such a lollipop diagram. Thus, if $\Area(\Delta)$ denotes the number of 2-cells in $\Delta$ and $\Area(w)$ denotes the minimum of $\Area(\Delta)$ over all van~Kampen diagrams $\Delta$ for $w$, then $\Area(w)$ also equals the minimum $\ell$ over all free equalities of the form \eqref{prod}.  (For a careful treatment, see \cite{Bridson6}.)

Similarly, \emph{diameter} has both a diagrammatic definition and an algebraic one. (We consider the \emph{intrinsic} diameter here, rather than its \emph{extrinsic} counterpart, in which distances are measured in the Cayley graph \cite{BR1}.) For a van~Kampen diagram $\Delta$, let $\Diam(\Delta)$ denote the minimum $L$ such that for all vertices $v$ in $\Delta$, there is a path of at most $L$ edges in $\Delta$ connecting $v$ to the base vertex $\ast$.
Then $\Diam(w)$ is the minimum of $\Diam(\Delta)$ among all van~Kampen diagrams $\Delta$ for $w$. Define $\Diam'(w)$ to be the minimum $L$ such that there is a free equality of the form \eqref{prod} with $|u_i| \leq L$ for all $i$. Then
$$\Diam'(w) \leq \Diam(w) \leq \Diam'(w)+K+|w|,$$
where $K = \max_{r \in R}|r|$. 
To see this, let $\Delta$ be a van Kampen diagram for $w$ with diameter equal to $\Diam (w)$.  Then there exists a lollipop diagram $\Delta'$, obtained by cutting open $\Delta$ as described above, so that the map $\Delta' \to \Delta$ sends the edge-paths $\rho_i$ labelled by the $u_i$ to geodesic paths in the 1-skeleton of $\Delta$ ending at $\ast$.  Then the first inequality is immediate.  
For the second inequality, let $\Delta$ be a van Kampen diagram for $w$ that is obtained by folding up a lollipop diagram arising from an expression \eqref{prod}  realizing $\Diam'(w)$.  Then every vertex $v$ of $\Delta$ lies either on the boundary of $\Delta$ or on the boundary of a face of $\Delta$. In the first case, a portion of the boundary circuit gives a path from $v$ to $\ast$ of length at most $|w|$ (in fact, at most $|w|/2$). In the second case, such a path can first follow a portion of the boundary circuit of the face and then can follow the image in $\Delta$ of one of the paths $\rho_i$.

\emph{Radius} is a variant of diameter: we define $\Rad(\Delta)$ to be the minimum $L$ such that every vertex $v$ in $\Delta$ can be connected to the boundary of $\Delta$ by a path of at most $L$ edges. If the boundary circuit of $\Delta$ has length $n$, then
\begin{equation} \label{eq diam rad}
 \Rad(\Delta) \leq \Diam(\Delta) \leq \Rad(\Delta) + \lfloor n/2 \rfloor, 
 \end{equation} 
with the first inequality holding because $\ast$ is on the boundary, and the second because an edge-path to the boundary can be extended to an edge-path to $\ast$ by traversing a portion of the boundary circuit in one direction or the other.  As before, $\Rad(w)$ is defined to be the minimum of $\Rad(\Delta)$ among all van~Kampen diagrams $\Delta$ for $w$.

The \emph{Dehn function} (or \emph{minimal isoperimetric function}) $\Area :  \N \to \N$, the \emph{diameter function} (or \emph{minimal isodiametric function}) $\Diam : \N \to \N$, and the \emph{radius function} $\Rad : \N \to \N$   are then defined so that $\Area(n)$,  $\Diam(n)$, and $\Rad(n)$ are, respectively, the maximum of $\Area(w)$, $\Diam(w)$,  and $\Rad(w)$,   over all words $w$ of length at most $n$ that represent $1$ in $G$. (We write $\Area_G$ or similar when we wish to resolve ambiguity over which group or which presentation we refer to.)

It follows from \eqref{eq diam rad} that for all $n \in \N$,
\begin{equation} \label{eq diam rad functions}
 \Rad(n) \leq \Diam(n) \leq \Rad(n) + \lfloor n/2 \rfloor.
\end{equation}

A van~Kampen diagram $\Delta$ is \emph{reduced} when  it does not contain a pair of   back-to-back cancelling 2-cells, which is to say that there is no edge $e$ in $\Delta$ which is in the boundary of two 2-cells in such a way that the same word is read around those 2-cells starting from $e$, one clockwise  and the other anticlockwise.   If a word represents the identity, then it admits a reduced van~Kampen diagram.

\subsection{Area-diameter pairs} \label{sec: ad pairs}

A pair $(f,g)$ of non-decreasing functions $\N \to \N$ is an \emph{area-diameter pair for $G$} if, for every $n \in \N$ and every word $w$ of length $|w| \leq n$ that represents $1$, there is a van~Kampen diagram $\Delta$ with $\Area(\Delta) \leq f(n)$ and $\Diam(\Delta) \leq g(n)$. 
It follows that $w$ can be expressed as in \eqref{prod} with

\begin{equation}  \label{bounds by f and g}
\ell  \leq f(n)  \qquad  \text{ and }   \qquad     \max \{ |u_i| \mid 1 \leq i \leq {\ell} \}  \leq g(n).
\end{equation} 

\begin{remark} \label{rem: diam}
Suppose $(f, g)$ is an area-diameter pair for $G$.  Then for 
all $n \in \N$,
$$g(n) \geq \Diam(n) \geq \lfloor n/2 \rfloor,$$
 assuming, for the second inequality, that $G$ is infinite.  
 The first inequality is an easy consequence of the definitions of $\Diam(n)$ and $g(n)$.  
 The second holds because for all $n$, there exists an element $b_n \in G$ such that $d_S(1, b_n) = \lfloor n/2 \rfloor$.  Let $u_n$ be a minimal length word representing $b_n$. Suppose $\Delta$ is any van~Kampen diagram for $w_n = u_n u^{-1}_n$. 
 Then around $\partial \Delta$, starting from  the base vertex $\ast$, we read $u_n$ and then $u^{-1}_n$, and   the distance in the 1-skeleton of $\Delta$ from $\ast$ to the vertex  where we switch from  $u_n$ to $u_n^{-1}$ must be at least $|u_n|$. So  $\Diam(\Delta) \geq \lfloor n/2 \rfloor$, hence $\Diam(n) \geq \lfloor n/2 \rfloor$.
\end{remark}

\subsection{\texorpdfstring{$t$-corridors}{t-corridors}} \label{sec: corridors}

Per the HNN structure we will establish in Corollary~\ref{cor:hnn structure}, the presentation \eqref{Rips presentation} for $G$ can be expressed as  $\langle \, \S \mid t U_i t^{-1} = V_i, 1 \leq i \leq k + 6m \rangle$, where $U_i$ and $V_i$ are words that contain no $t$-letters.  So every 2-cell   in a van~Kampen diagram  $\Delta$ has exactly two edges labelled $t$, and  the 2-cells in $\Delta$ concatenate along common $t$-edges  to make \emph{$t$-corridors} (when they terminate at a pair of $t$-edges in $\partial \Delta$) or \emph{$t$-annuli} (when they close up).  Along each side of a $t$-corridor there is an edge-path, and  
 around each side of a $t$-annulus there is an edge-loop.  Assuming those edge-paths and edge-loops  start from suitable basepoints, on one side they are labelled  by a concatenation of the $U_i^{\pm 1}$ and on the other by a concatenation of the $V_i^{\pm 1}$.  They need not be injective paths. The words read along them need not be reduced.

\section{Small cancellation} \label{sec: small cancellation}

Recall that $\X$ denotes the set of the $2k+12m$ variously decorated words $X$ that appear in the presentation $\langle \S \mid  \R\rangle$ for $G$ given in  \eqref{Rips presentation} and these words are disjoint subwords of the infinite Rips word $x_1 x_2^1 \, x_1 x_2^2 \, x_1 x_2^3 \cdots$. Here we discuss how these subwords should be chosen and record some consequences.

Define 
\begin{equation} \label{Ks}
\lambda = \max \set{|r_i| \mid 1 \leq i \leq k }  \  \text{ and } \ \lambda' =  \max \set{\lambda+2, 5}.
\end{equation}
Building on \cite[Remark~3.2]{Wise1}, we define $\xi$ to be any positive integer such that
\begin{equation} \label{ineq:small} 
  \alpha = \xi (2k+12m)+\xi-2   \ \text{ satisfies } \      
4 \alpha  +\lambda'  < (3 \xi +1) \xi /8.
\end{equation}
Such a $\xi$ exists because 
 $\alpha$ grows linearly as a function of $\xi$, while $(3 \xi +1) \xi /8$ grows quadratically.
Then we take 
\begin{equation} \label{eq:Chi}
  \X = \set{\chi_i \mid 1 \leq i \leq 2k+12m}, \qquad \textup{where} \qquad \chi_i = x_1 x_2^{\xi i}x_1 x_2^{\xi i+1} \cdots x_1 x_2^{\xi i+\xi-1}.
  \end{equation}
It is immaterial which decorated word $X$ is assigned to each $\chi_i$.

A set $\mathcal{B}$ of words satisfies the \emph{uniform $C'(\zeta)$ condition} if every common prefix (a \emph{piece}) of a pair of distinct cyclic permutations of words in $\mathcal{B}^{\pm 1}$ has length strictly less than $\zeta$ times the length of the shortest word in $\mathcal{B}$.

\begin{lemma} \label{lem:C'1/4}
Let $\X$ be the set defined in~\eqref{eq:Chi}, let $\U$ and $\V$ be the sets of words defined in~\eqref{eq:UV}, and let $\R$ be the set of defining relators of $G$ in \eqref{Rips presentation}.  Then
\begin{enumerate}
  \item \label{X} $\X$ satisfies the uniform $C'(1/4)$ condition.
  \item \label{UV} $\U \cup \V$   satisfies the uniform $C'(1/4)$ condition. 
  \item \label{R} $\R$ satisfies the uniform $C'(1/6)$ condition.
\end{enumerate}
\end{lemma}

\begin{proof} The word $\Omega = x_2^{\alpha-1} x_1 x_2^{\alpha}$ is a longest piece for $\X$.  (This uses that $\xi \geq 3$, which \eqref{ineq:small}  forces.)
It  occurs twice as a subword of $\chi_{2k+12m}$ and has length $2\alpha$.
The shortest word in $\X$ is $\chi_1$, which has length $|\chi_1| = (3 \xi +1)\xi/2$. We have 
\begin{equation} \label{eq:2alpha}
2\alpha \leq 4\alpha + \lambda' < |\chi_1|/4
\end{equation}
 by \eqref{ineq:small}, and so $\X$ satisfies the uniform $C'(1/4)$ condition, proving \eqref{X}.

The constant $\lambda'$ of \eqref{Ks} is defined so that every relation of \eqref{Rips presentation} contains at most $\lambda'$ letters (counted with repetitions) outside the subwords denoted by variously decorated $X$ symbols.  A piece for $\U \cup \V$ or  $\R$ must be a concatenation of two pieces  from $\X$ and words made up of some of those $\lambda'$ letters, so the longest  piece has length at most  $2|\Omega| + \lambda'= 4 \alpha + \lambda'$.

The shortest word in $\U \cup \V$ has length at least $|\chi_1|$, so \eqref{UV} follows from \eqref{eq:2alpha}. 
The shortest defining relator in $\R$ has length at least $2|\chi_1| = (3\xi+1)\xi$.  For \eqref{R}, again using \eqref{eq:2alpha}, we have $4\alpha + \lambda'< |\chi_1|/4  < 2|\chi_1|/6$, as required.
\end{proof}

Therefore the presentation $G = \langle \S \mid  \R \rangle$ is $C'(1/6)$, so $G$ has a linear isoperimetric inequality and is hyperbolic \cite{Gersten9}. In fact, the \emph{uniform} $C'(1/6)$ condition gives more: $G$ is $\textup{CAT}(-1)$, cubical, virtually special, and residually finite \cite{Agol, BrownS, Gromov7, Martin, WiseCubulating}.

Our next lemma is a consequence of the $C'(1/4)$ condition, and immediately implies  the corollary that  follows.

\begin{lemma} \label{C14 applications}
Let $\U$ and $\V$ be the sets of words 
defined in~\eqref{eq:UV}.  
Then both $\U$ and $\V$  freely generate free subgroups of  $F = F(a_1, \ldots, a_m, x_1, x_2)$. 
\end{lemma}

\begin{proof}
Suppose $w = Y_1 \cdots Y_l$ is any non-empty freely reduced word on $\U^{\pm 1}$ or on $\V^{\pm 1}$. To prove the result, we must show that $w$, viewed 
as a word on $a_1, \ldots, a_m, x_1, x_2$, cannot freely reduce to the empty word. Since $\U$ and $\V$ satisfy  the uniform $C'(1/4)$ condition by Lemma~\ref{lem:C'1/4}(2),  free reduction between $Y_i$ and $Y_{i+1}$ erodes fewer than a quarter of the letters of each.  Thus, for each $i$, at least half the letters of $Y_i$ survive cancellation with its neighbors.  Therefore $w$ cannot freely equal the empty word.
\end{proof}

\begin{corollary}\label{cor:hnn structure}
  The group $$G = \langle \S \mid t U_i t^{-1} = V_i, 1 \leq i \leq k + 6m \rangle$$ is an HNN extension of the free group $F = F(a_1, \ldots, a_m, x_1, x_2)$ with stable letter $t$ and associated subgroups generated by the free bases $\U = \set{U_i}$ and $\V = \set{V_i}$.   
\end{corollary}

\section{The upper bound on distortion} \label{sec: upper bound}

In this section we prove the upper bound in Theorem~\ref{Main}.  Let $Q = \langle a_1, \dots, a_m \mid r_1, \dots, r_k \rangle$, $G=\langle \S \mid \R \rangle$ and $N$ be the groups from Section~\ref{sec: Wise HNN},  
and let $(f,g)$ be an area-diameter pair for $Q$.  Let $K= \max \set{ |r| \mid r \in \R }$, the length of the longest defining relator in our presentation for $G$. 

Suppose $w$  is a word of length at most $n$ on $\S$ that represents an element of $N \leq G$.  We will produce a word on $\{ t, x_1, x_2 \}$ representing $w$, whose length is bounded above by a suitable  function  
involving $f$ and $g$, 
leading to the desired upper bound on the distortion of $N$.

Let $\wo$ be $w$ with all letters  $t, x_1, x_2$ deleted.  Because $N$ is normal in $G$, the $|\wo|$ instances of the letters $a_1, \ldots, a_m$ in $w$ can be shuffled to the right-hand end of $w$ without changing the group element represented.  Initially, there are $(|w| - |\wo|)$ letters  $t, x_1, x_2$ in $w$, and for all $i$, the words $a_i^{\pm 1} t a_i^{\mp 1}$, $a_i^{\pm 1} x_1 a_i^{\mp 1}$, and $a_i^{\pm 1} x_2 a_i^{\mp 1}$ equate in $G$ to words of length at most $K$ on $\{t, x_1, x_2\}$.  So  there is a word $Y$ on $\{t, x_1, x_2\}$ of length   
\begin{equation} \label{Y length bound}
|Y| \leq (|w| - |\wo|) K^{|\wo|} \leq n K^n
\end{equation} 
 such that 
\begin{equation} \label{Ywow}
w = Y \wo  \ \text{ in } G.  
\end{equation} 

Moreover, $\wo = 1$ in $Q$, because $w$ represents an element of $N$, and so there is a word 
\begin{equation} \label{W}
W    =    \left( {u_1}^{-1} r_{j_1}^{\epsilon_1} u_1 \right) \cdots  \left( {u_{\ell}}^{-1} r_{j_{\ell}}^{\epsilon_{\ell}} u_{\ell} \right)
\end{equation} 
such that    $j_1, \ldots, j_{\ell} \in \{1, \ldots, k \}$,    $\epsilon_1, \ldots, \epsilon_{\ell} \in \{1, -1\}$, each $u_i$ is a word on $\{a_1, \dots, a_m\}$,
\begin{equation} \label{Wwo}
W = \wo \ \text{ in } F(a_1, \ldots, a_m),
\end{equation} 
and the bounds  \eqref{bounds by f and g} hold, i.e.~$\ell  \leq f(n)$ and $\max \{ |u_i| \mid 1 \leq i \leq {\ell}\} \le g(n)$.

Viewing $W$ as a word on $\S$, we will use defining relations from $\mathcal R$ to express it in terms of $t$, $x_1$, and $x_2$. 
From the first family of relations in \eqref{Rips presentation}, we see that for every $i$, we have $r_{j_i} = X_{j_i} t X'_{j_i} t^{-1}$ in $G$.  Observe that 
the word on the right has length at most $K$. If we use the second and third families of defining relations in \eqref{Rips presentation} to conjugate $X_{j_i} t X'_{j_i} t^{-1}$ repeatedly by letters in $a^{\pm 1}_1, \ldots, a^{\pm 1}_m$, its length as a word on $t, x_1, x_2$ increases by at most a factor of $K$ at each step. Since $|u_i| \leq g(n)$ and $|X_{j_i} t X'_{j_i} t^{-1}| \le K$, it follows that ${u_i}^{-1} r_{j_i}^{\epsilon_i} u_i$ is equal, in $G$, to a word $Z_i$ on $\{t, x_1, x_2\}$ of length
\begin{equation} \label{length of Zi}
|Z_i| \leq K^{g(n)+1}.  
\end{equation} 

Combining this with \eqref{Ywow} and \eqref{Wwo}, we learn that
\begin{equation} \label{w in N}
w = Y \wo =  Y W  =  Y Z_1 \cdots Z_{\ell}  \ \text{ in } G.
\end{equation} 
This equates $w$ with a word $Y Z_1 \cdots Z_{\ell}$ on the generators $t, x_1, x_2$ of $N$. By \eqref{Y length bound}, \eqref{length of Zi}, and $\ell \leq f(n)$, 
$$|Y Z_1 \cdots Z_{\ell}|   \leq n K^{n} + f(n) K^{g(n)+1}.$$

Since Theorem~\ref{Main} assumes that $Q$ is infinite, Remark~\ref{rem: diam} applies, and implies that $n \leq g(2n)$ for all $n \in \N$.
 Then $K^n \leq K^{g(2n)}$, and so $K^n \preceq K^{g(n)}$.
It follows that 
$$n K^{n} + f(n) K^{g(n)+1} \preceq n K^{g(n)} + f(n) K^{g(n)+1} = (n + K f(n)) K^{g(n)}.   
$$
Thus $\Dist_N^G(n) \preceq (n+f(n)) \,K^{g(n)}$, as required.

 \section{The lower bound on distortion}  \label{sec: lower bound}
In this section we will prove the lower bound in Theorem~\ref{Main}. 
Let $1 \to N \to G \to Q \to 1$ be as in the statement of the theorem. 
We will show that the distortion of $N$ in $G = \langle \S \mid \R \rangle$ satisfies
\begin{equation}\label{eq:lower bound}
K_1^{\left(\Diam_Q(n) - \lfloor n/2 \rfloor\right)} \preceq \Dist_N^G(n)
\end{equation}
for a suitable constant $K_1 > 1$.
This will be a consequence of the following proposition, concerning the radius function. Define 
\begin{equation} \label{eq: constants}
C = 2 \max \set{ |w| \mid w \in \U \cup \V }, \qquad  K= \max \set{ |r| \mid r \in \R }, \qquad M = C+K. \end{equation} 
Note that $M > 0$, because $Q$ is a non-trivial group (since it is assumed to be infinite), implying \eqref{Rips presentation} has at least one defining relation.

\begin{prop} \label{prop lower bound}  
If $w$ is a word on $\S$ 
that equals a word $v$ on $\{t, x_1, x_2\}$ in $G$, and if $\overline w$ is the word obtained from $w$ by deleting all instances of $t$, $x_1$, and $x_2$,  then 
\begin{equation}\label{v lower bound}
|v| \ge 2^{\frac{1}{M} \Rad_Q (\overline w)}  -|w| -2. 
\end{equation}

\end{prop}

Before proving Proposition~\ref{prop lower bound}, we explain how it yields the desired lower bound in Theorem~\ref{Main}:
\begin{proof}[Proof of the lower bound \eqref{eq:lower bound}]
For each $n \geq 2$, we will produce a word $w$ on $\S$ of length at most $n$ representing an element of $N$, such that every word on $\{t, x_1, x_2\}$ representing $w$ has length bounded below by a function
$\simeq K_1^{\left(\Diam_Q(n) - \lfloor n/2 \rfloor\right)}$, for a suitable $K_1>1$.  
To that end, for each $n\ge 2$ choose a word $\overline w$ on $\{a_1, \dots, a_m\}$ that represents the identity in $Q$, has length at most $n$, and realizes the diameter function of $Q$, meaning that every van~Kampen diagram for $\overline w$ over $Q$ has diameter at least $\Diam_Q(n)$. By
~\eqref{eq diam rad}, $ \Rad_Q (\overline w) \ge \Diam_Q(n) - \lfloor n/2 \rfloor $.

The word $w:=\overline w$, regarded as a word in the generators of $G$, maps to $\overline w$ under the quotient map from $G$ to $Q$. Since $\overline w$ represents $1$ in $Q$, the word $w$ represents an element of $N$. Let $v$ be any word on $\{t, x_1, x_2\}$ that equals $w$ in $G$. Then, 
taking $K_1 = 2^{\frac 1 M}$,
Proposition~\ref{prop lower bound} gives
$$|v| \ge 2^{\frac{1}{M}\Rad_Q (\overline w)}  -|w| -2
\ge  2^{\frac{1}{M}\left(\Diam_Q(n) - \lfloor n/2 \rfloor\right)} -n -2
\succeq K_1^{\left(\Diam_Q(n) - \lfloor n/2 \rfloor\right)}. 
$$
This proves the lower bound~\eqref{eq:lower bound}. (Both $\Dist_N^G(1)$ and $\Diam_Q(1)$ are constants, so we may ignore $n \leq 1$.)
\end{proof}

We work toward our proof of Proposition~\ref{prop lower bound} through four lemmas.
The first allows us to ``collapse'' a van~Kampen diagram  over  $G = \langle \S \mid \R \rangle$ to one over $Q = \langle a_1, \dots, a_m \mid r_1, \dots, r_k \rangle$. 

\begin{lemma} \label{diagram collapse}
Let $u$ be a word on $\S$ that represents the identity in $G$, and let
$\Delta$ be a van~Kampen diagram for $u$ over $G$. Let $\overline u$ be the word obtained from $u$ by deleting all instances of $t$, $x_1$, and $x_2$. Then there exist a van~Kampen diagram $\overline \Delta$ for $\overline u$ over $Q$ and a quotient map $\Delta \to \overline \Delta$ such that $\partial \Delta$ maps to $\partial \overline \Delta$ and: 

\begin{enumerate}
\item  \label{property 1}
Every 2-cell of $\overline \Delta$ is the image of a 2-cell of $\Delta$.
\item \label{property 2} Composing with the quotient map $\Delta \to \overline \Delta$ takes every edge-path $\gamma$ in the 1-skeleton of $\Delta$  to an edge-path in the 1-skeleton of $\overline \Delta$ of combinatorial length (meaning number of edges it traverses) less than or equal to that of $\gamma$.  
\end{enumerate}
\end{lemma}

\begin{proof}
We inductively construct a sequence of quotient maps
$$
\Delta = \Delta_0 \to \Delta_1 \to \Delta_2 \to  \Delta_3 \to  \dots
$$
such that each $\Delta_i$ is a planar 2-complex, has each edge labelled by a generator of $G$, and $\Delta \to \Delta_i$ enjoys properties \eqref{property 1} and \eqref{property 2}.  
For $i>0$, the complex $\Delta_i$ is obtained from $\Delta_{i-1}$ via one of the following moves:

\begin{enumerate}
\item[(i)] \label{} Collapse an edge labelled $t$, $x_1$, or $x_2$ whose endpoints are distinct vertices; 
\item[(ii)] Collapse a subcomplex $D \subset \Delta_{i-1}$ to a single point, provided that $D$ is homeomorphic to a 2-disc and its boundary is a single edge (necessarily with its endpoints identified) labelled $t$, $x_1$, or $x_2$.
\end{enumerate}

Because of the HNN structure of $G$ (Corollary~\ref{cor:hnn structure}), the elements $t$, $x_1$, and  $x_2$ have infinite order, so no edge in $\Delta$  labelled $t$, $x_1$, or $x_2$ is a loop. Thus $\Delta_1$ must be obtained from $\Delta_0$ by a move of type (i). A move of type (ii) is necessary at some stage if $\Delta$ contains a subdiagram whose boundary is labelled by a word on $\{t, x_1, x_2\}$.

Each $\Delta_i \to \Delta_{i+1}$ is the collapse of a contractible subcomplex. The nature of the moves (i) and (ii) guarantees that $\Delta_{i+1}$ is a  planar 2-complex because $\Delta_i$ is planar;  the edges in $\Delta_{i+1}$ inherit labellings by  letters in $\S$ from such labellings in $\Delta_i$. Each 2-cell of $\Delta_{i+1}$ is the image of a 2-cell $c$ of $\Delta_i$ whose boundary label is either unchanged  or has one letter deleted.

The process stops at $ \Delta_l$ for some $l$ when there are no remaining edges labelled $t$, $x_1$, or $x_2$. The boundary word of each of its 2-cells agrees with that of some 2-cell of $\Delta$ after all instances of $t$, $x_1$, and $x_2$ have been deleted.  Thus, after cyclic permutation and inversion as necessary, such a boundary label is either $r_i$ for some $1 \le i \le k$ or $a_ja_j^{-1}$ for some $1 \le j \le m$.  If there is a cell of the latter type, we collapse it to a single oriented edge labelled $a_j$.  This can be done by identifying the two edges labelled $a_j$ in its boundary by an orientation-preserving homeomorphism, and then extending continuously over the interior of the cell. 
The result of this collapse is again a planar 2-complex.  We call this a type (iii) move.  

After finitely many type (iii) moves, we end up with a planar 2-complex $\overline \Delta$ with oriented edges labelled by generators of $Q$, such that each of its 2-cells has boundary labelled by a relator of $Q$.  
Likewise, the word read around $\partial \overline \Delta$ is obtained by deleting all instances of $t$, $x_1$, and $x_2$ from the word read around $\partial \Delta$. 
Thus $\overline \Delta$ is a van~Kampen diagram for $\overline u$ over $Q$, and the quotient map $\Delta \to \overline \Delta$ satisfies \eqref{property 1} and takes $\partial \Delta$ to $\partial \overline \Delta$.

Composing with the successive quotient maps takes an edge-path $\gamma$ in $\Delta$ to a sequence of paths $\gamma_i$ in $\Delta_i$.   If $\Delta_{i-1} \to \Delta_i$ is a type (i) move, then the word read along $\gamma_i$ is that read along $\gamma_{i-1}$ but with letters $t$,  $x_1$, or  $x_2$  removed, where $\gamma_{i-1}$ runs along the edge that is to be collapsed.    If $\Delta_{i-1} \to \Delta_i$ is a type (ii) move, then the word read along $\gamma_i$ is that read along $\gamma_{i-1}$ but with some subwords removed --- more specifically, the subwords where $\gamma_{i-1}$ enters the subcomplex that is collapsed.  (These subwords could feature any of the letters in $\S^{\pm 1}$.)  If $\Delta_{i-1} \to \Delta_i$ is a type (iii) move, then the label is unchanged. So 
\eqref{property 2} holds.
\end{proof}

\begin{lemma}\label{t corridors}
Let $\Delta$ be a van Kampen diagram over $G$.  Then the following hold.
\begin{enumerate}
\item Every 2-cell of $\Delta$ is contained in a $t$-corridor or $t$-annulus.
\item If $\Delta$ is reduced, then there are no $t$-annuli.
\end{enumerate}
\end{lemma}

\begin{proof}
The first property is immediate. The second follows from the HNN structure for $G$ described in Corollary~\ref{cor:hnn structure}, as we now explain.

If a $t$-annulus exists, consider an innermost one. The inner boundary of this  annulus encloses no $2$-cells; such a $2$-cell would have to be contained in a $t$-corridor or $t$-annulus by (1), which is impossible since the inner boundary has no edges labelled $t$, and the chosen annulus was an innermost one. 
The word read around it is therefore a word on $\U$ or on $\V$, depending on the orientation of the $t$-edges in the annulus. Moreover, it is freely reducible as a word on $\U$ or $\V$, by Lemma~\ref{C14 applications}, so it has a subword of the form $U_i^{\pm 1} U_i^{\mp 1}$ or $V_i^{\pm 1} V_i^{\mp 1}$ for some $U_i \in \U$ or $V_i \in \V$. This subword is read along portions of the boundaries of two cells in the annulus that share a $t$-edge. Since each $U_i$ and $V_i$ appears in exactly one defining relation in the presentation $\langle \S \mid \R \rangle$ for $G$, these adjacent 2-cells form a cancelling pair of the same type but opposite orientations, contradicting the hypothesis that $\Delta$ is reduced.
\end{proof}

Recall that $C = 2 \max \set{ |w| \mid w \in \U \cup \V }$.  
Our next lemma says that two $t$-corridors in a  reduced van~Kampen diagram over $G=\langle \S \mid \R \rangle$ cannot travel side-by-side for a distance $C$ or more, and a single $t$-corridor cannot loop back on itself beyond a distance of $C/2$.  To state this precisely, we use the terminology set out in Section~\ref{sec: corridors}.

\begin{lemma}\label{key corridor lemma}  
  If $\Delta$ is a reduced van Kampen diagram over $G$, then the following hold:
  \begin{enumerate}
  \item \label{corridors part 1}
  Suppose $\tau_1$ and $\tau_2$ are distinct $t$-corridors and that $u$ and $v$ are vertices in $\tau_1 \cap \tau_2$.  For $i=1,2$, let $\rho_i$ be the subpath from $u$ to $v$ of an edge-path along the appropriate side of $\tau_i$.   
    Then the lengths of $\rho_1$ and $\rho_2$ are both strictly less than $C$.
  
  \item \label{corridors part 2} If a subpath $\rho$ of 
  an edge-path along a side of a $t$-corridor is a loop, then $|\rho| < C/2$. 
  \end{enumerate}
\end{lemma}

\begin{proof}
The labels of the boundaries of $t$-corridors are concatenations of elements (we call these \emph{syllables}) of $\U^{\pm 1}$ or $\V^{\pm 1}$.  So  the labels of $\rho_1$ and $\rho_2$ are words on $\mathcal \U$ or on $\mathcal V$, possibly with the initial or final syllable truncated.  Accordingly, we denote the label of $\rho_1$ by $\P = \widehat{\P}_1 \P_2 \cdots \P_{p-1} \widehat{\P}_p$ and that of $\rho_2$ by $\Q = \widehat{\Q}_1 \Q_2 \cdots \Q_{q-1} \widehat{\Q}_{q}$, where $\{\P_1, \dots, \P_p\}$ and $\{\Q_1, \dots, \Q_q\}$ are subsets of either $\U^{\pm 1}$ or $\V^{\pm 1}$, and $\,\widehat{Z}\,$ denotes a non-empty suffix or prefix of a syllable $Z$ as appropriate.

The circuit $\omega = \rho^{-1}_1 \rho_2$ encloses no 2-cells. Indeed, any 2-cell inside $\omega$ would have to belong to a $t$-annulus --- it could not belong to a $t$-corridor because $\omega$ has no $t$-edges --- but Lemma~\ref{t corridors} rules out $t$-annuli in a reduced diagram.  It follows that $\rho_1$ is homotopic to $\rho_2$ in the 1-skeleton $\Delta^{(1)}$ of $\Delta$.  More concretely, there is an embedded edge-path $\mu$, 
 a homotopy from $\rho_1$ to $\mu$ in $\Delta^{(1)}$ which corresponds to free reduction among successive syllables of $\P$, and a homotopy from $\rho_2$ to $\mu$ in $\Delta^{(1)}$ which corresponds to free reduction among successive syllables of $\Q$. Thus $\mu$ inherits a label 
$\P' = {\P}_1' \P_2' \cdots \P_{p-1}' {\P}_p'$ from the free reduction of $\P$ and  $\Q' = {\Q}_1' \Q_2' \cdots \Q_{q-1}' {\Q}_{q}'$ from the free reduction of $\Q$, and these two labels are identical as words.

By Lemma~\ref{lem:C'1/4}(2), $\U \cup \V$ satisfies the uniform $C'(1/4)$ condition.  Let $L$ denote the length of the longest piece of $\U \cup \V$.  Then the length of each $\P_i$ or $\Q_i$ is strictly greater than $4L$.  Moreover, for each $i$, we have $\P_i \neq \P_{i+1}^{-1}$ and $\Q_i \neq \Q_{i+1}^{-1}$, for otherwise $\Delta$ would not be reduced. 
So at most $2L$ letters of $\P_i$ or $\Q_i$ are cancelled away during the free reduction procedure and we have that 
\begin{equation} \label{PiQi ineq}
|\P_i'|> 2L  \text{ for all } 2 \leq i \leq p-1 \text{ and } |\Q_i'|>2L  \text{ for all } 2 \leq i \leq q-1. 
\end{equation}

\medskip
Now suppose that one of the paths, say $\rho_1$, has length at least $C$.  Since $|\P_i| \leq  C/2$ for all $i$, either $p \ge 3$ or $\P = \P_1\P_2$ with neither syllable truncated (which is to say, $\widehat{\P_1} = \P_1$ and $\widehat{\P_2} = \P_2$). In either case, $\P_2$ is  a subword of $\P$. 

\medskip
\noindent \emph{Claim.}  
Given that $\mathcal P_2$ is a subword of $\mathcal P$, the labels $\P'$ and $\Q'$ of $\mu$ agree on a subpath $\mu'$ of length strictly greater than $L$, such that $\mu'$ is in the boundary of a single cell $c_1$ of $\tau_1$ and a single cell $c_2$  of $\tau_2$.

\medskip

To prove the claim, we first observe that  $|\P_2'|> 2L$, because if $p \geq 3$ this is an instance of \eqref{PiQi ineq}, and if $\P = \P_1 \P_2$ then no more than $L$ letters of $\P_2$ can be cancelled away (with $\P_1$ necessarily) on free reduction of $\P$. 
 
Since $\P'$ and $\Q'$ are  identical as words, $\P_2'$ is a subword of $\Q'$, and so there exist $j \ge 1$ and $l\ge 0$ such that $\P_2' = \widehat{\Q_j'} \Q_{j+1}' \cdots \widehat {\Q_{j+l}'}$ as words, where $\widehat{\;}$ denotes possible truncation as before. We examine three cases:

\medskip

\noindent \emph{Case $l=0$.} We have that $\P_2' = \widehat{\Q_j'}$ and so $\P'_2$ is a subword of $\Q_j$, and we take $\mu'$ to be the subpath of $\mu$ labelled by this $\P_2'$.

\noindent \emph{Case $l=1$.} At least one of $ \widehat{\Q_j'}$ and $\widehat{\Q_{j+1}'}$ agrees with half or more of $\P_2'$, and we take $\mu'$ to be the subpath of $\mu$ of that longer overlap. 

\noindent \emph{Case $l \geq 2$.} In this case $\Q_{j+1}'$ is a subword of $\P_2'$, and we take  $\mu'$ to be the subpath of $\mu$ labelled by $\Q_{j+1}'$.

\medskip

\noindent That the length of $\mu'$ is strictly greater than $L$ follows from  $|\P_2'|> 2L$ in the first two cases and from $|\Q_{j+1}'| > 2L$ by   \eqref{PiQi ineq} in the third case.

\medskip

 Having established the claim, we can complete the proof of Part~\eqref{corridors part 1}: were the lengths of $\rho_1$ and $\rho_2$ not both strictly less than $C$, then the claim would tell us that $\partial c_1$ and $\partial c_2$ agree on a  subpath of length longer than that of any piece, so they would carry the same defining relation but with opposite orientations, and therefore contradict the assumption  that $\Delta$ is reduced.   
 
\medskip

 For Part~\eqref{corridors part 2}, suppose the label of $\rho$ is $\widehat{\P}_1 \P_2 \cdots \P_{p-1} \widehat{\P}_p$, with $\widehat{\P_1}, \widehat{\P_p}$ non-trivial, and (as $\Delta$ is a reduced diagram) $\P_i \neq \P_{i+1}^{-1}$  for each $i$.  We must have $p\ge 2$: elements of $\mathcal U$ and $\mathcal V$ are reduced words in $F(a_1, \dots, a_m, x_1, x_2)$, which embeds in $G$ by Corollary~\ref{cor:hnn structure}, so subwords of such elements cannot label loops. 
If $p=2$, then $\widehat{\P}_1\widehat{\P}_2$ freely reduces to the empty word, which can only happen if $\widehat{\P}_1= \widehat{\P}_2^{-1}$.  However, as ${\P}_1$ and ${\P}_2^{-1}$ are distinct, it follows that their overlap is at most $L$.  Thus the length of $\rho $ is at most 
$2L < C/2$.  Now if $p>2$, then $\P_2$ is untruncated.  
   We view $\rho$ as two non-trivial paths between the same pair of points, $\rho_1$ with label $\P_2 \cdots \P_{p-1}$ and $\rho_2$ with label $\widehat \P_1^{-1} \widehat{\P}_p^{-1}$. 
  As before $\rho_1$ and $\rho_2$ are homotopic in $\Delta^{(1)}$ to the same path which inherits the (freely reduced) label  $\P_2' \cdots \P_{p-1}'$ from $\rho_1$ and 
  ${\widehat{\P}'_1}{}^{-1}{\widehat{\P}'_p}{}^{-1}$ from $\rho_2$, with $|\P_i'| > 2L$ for $2 \le i \le p-1$. Then ${\P_2}'$ agrees with a subword of either ${\widehat{\P}'_1}{}^{-1}$ or ${\widehat{\P}'_p}{}^{-1}$ of length strictly greater than $L$, and as before, this leads to a contradiction to the assumption
that $\Delta$ is reduced.    
\end{proof}

\begin{figure}[ht]
\begin{overpic}
{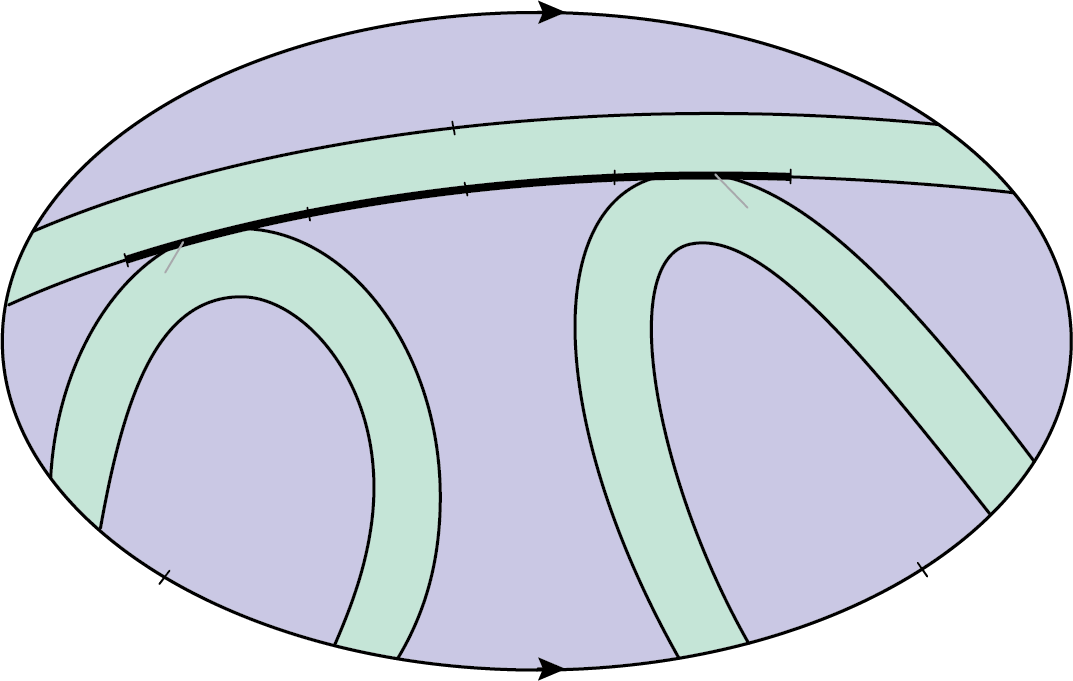}
\put(49.5,64.5){\small{$w$}}  
\put(49.5,-3){\small{$v$}}  
\put(10,58){\small{$\Delta$}}  
\put(41,54){\small{$p$}}  
\put(43,42){\small{$p'$}}  
\put(81,48){\small{$\tau$}}      
\put(35,24){\small{$\tau_0$}}      
\put(83,30){\small{$\tau_1$}}
\put(13.5,35.5){\small{$p_0$}}      
\put(70.5,42){\small{$p_1$}}      
\put(34.7,46.9){\small{$\mu'_0$}}   
\put(18,43.8){\small{$\mu''_0$}}   
\put(48,48.5){\small{$\mu'_1$}}   
\put(65,49){\small{$\mu''_1$}}   
\put(-1.5,38){\small{$t$}}  
\put(92.5,49){\small{$t$}}         
\put(5,13){\small{$t$}}         
\put(33,-1){\small{$t$}}  
\put(66,-1){\small{$t$}}         
\put(95,15){\small{$t$}}         
\end{overpic}
 \caption{Lemma~\ref{doubling corridors} illustrated}
  \label{fig:Delta}
\end{figure}

Our next lemma describes a situation in which the existence of additional $t$-corridors can be deduced from the existence of a single long one. It provides the inductive step in the proof of Proposition~\ref{prop lower bound}.  Its statement and proof are illustrated in Figure~\ref{fig:Delta}.

\begin{lemma}\label{doubling corridors}  
Let $\tau$ be a $t$-corridor in a reduced van~Kampen diagram $\Delta$ over $\langle \S \mid \R \rangle$.  Suppose there is a vertex $p$ of $\Delta$ on one side of $\tau$ such that $d(p, \partial \Delta) >M$, where $M$ is the constant defined in \eqref{eq: constants}.  Then there exist vertices $p_0$ and $p_1$ on the  side of $\tau$ not containing $p$, and $t$-corridors $\tau_0$ and $\tau_1$ containing $p_0$ and $p_1$, respectively, such that  
$\tau$, $\tau_0$, and $\tau_1$ are distinct corridors, and 
\begin{equation}
d(p, p_i) < M  \  \text{ for } \  i=0,1.
\end{equation}
\end{lemma}

\begin{proof}
Since $d(p, \partial \Delta) >M$, $p$ is not in $\partial \Delta$ and so there is a 2-cell of $\tau$ whose boundary contains $p$. Let $p'$ be a vertex on the boundary of this 2-cell that is on the other side of $\tau$.  Then $d(p, p')< K$, where $K$ is the constant defined in \eqref{eq: constants}.  Now, $d(p', \partial \Delta) >C$, because  $d(p, \partial \Delta) >M$ and $M=K+C$.  

Let $\mu$ denote the edge-path along the side of $\tau$ containing $p'$. As per Figure~\ref{fig:Delta}, for $i=0,1$, let $\mu_i'$ denote the length-$C/2$ subpaths of $\mu$ emanating from $p'$ in each direction, and let $\mu_i''$ be their length-$C/2$ continuations  further away from $p'$.  
Then the concatenations $\mu_i'\mu_i''$ both have length $C$ and, by our choice of $p'$, are disjoint from $\partial \Delta$. 
We claim that there is a vertex $p_i$ along $\mu_i''$ which is contained in the boundary of the closure of $\tau$.  If not, then $\mu_i''$ would be entirely encompassed within a fold between neighboring 2-cells of $\tau$, resulting in a loop of length $\ge C/2$ along a side of $\tau$, contradicting    
Lemma~\ref{key corridor lemma}(2). By construction, the subpath of $\mu$ from $p'$ to $p_i$ has length at least $C/2$, and  $p_i$ is not on $\partial \Delta$.  So
there is a $2$-cell $c_0$ that is not contained in $\tau$, whose boundary contains $p_i$.

The length of the subpath of $\mu$ from $p_0$ to $p_1$ is at least $C$, so they are distinct points by Lemma~\ref{key corridor lemma}(2). 
Moreover, 
$
d(p, p_i) \le d(p,p') + d(p', p_i) < K+C = M,
$ for $i=0, 1$,
as required. 

Now, by Lemma~\ref{t corridors}, there are 
 $t$-corridors $\tau_0$ and $\tau_1$ containing $c_0$ and $c_1$ (hence $p_0$ and $p_1$) respectively.  By construction,  $\tau_0 \neq \tau$  and $\tau_1 \neq \tau$.  It remains to show that $\tau_0 \neq \tau_1$. 
If $\tau_0 = \tau_1$, then $p_0$ and $p_1$ must lie on the same side of $\tau_0$; otherwise, one end of $\tau_0$ would be trapped between a (possibly empty) portion of $\tau$ and a portion of $\tau_0$ and could not reach the boundary.  
Then $\tau$ and $\tau_0$ are distinct corridors, and $p_0$ and  $p_1$ are distinct points in $\tau \cap \tau_0$ such that the edge-path along $\tau$ from $p_0$ to $p_1$ has length at least $C$, which contradicts Lemma~\ref{key corridor lemma}\eqref{corridors part 1}.
Thus $\tau_0 \neq \tau_1$. 
\end{proof}

We are now ready to complete the proof of Proposition~\ref{prop lower bound}. 

\begin{proof}[Proof of Proposition~\ref{prop lower bound}]
We have a word $w$ on $\S$ and a word $v$ on $t, x_1, x_2$ such that $w=v$ in $G$; let $\overline w$ denote the word obtained from $w$ by deleting all letters $t, x_1, x_2$. 
Let $\Delta$ be a reduced van~Kampen diagram for $w v^{-1}$ over $G$, let $\overline \Delta$ be the corresponding diagram for $\overline w$ over $Q$ constructed in Lemma~\ref{diagram collapse}, and let $R$ be the radius of $\overline \Delta$.

  We aim to show that 
\begin{equation} \label{eqn: the aim}
|v| + |w| \geq 2^{R/M} -2
\end{equation}
which will prove the proposition because, by definition,  $\Rad_Q(\overline{w}) \leq R$.  

We may assume that $R > M$, as \eqref{eqn: the aim} is immediate for $R \le M$.

Let $q$ be a vertex of $\overline \Delta$ such that a shortest path in the 1-skeleton from $q$ to $\partial \overline \Delta$ has length $R$, and let $p$ be a vertex of $\Delta$ mapping to $q$.  Then by Lemma~\ref{diagram collapse}(2),
\begin{equation} \label{pqR ineq} 
d(p, \partial \Delta) \geq d(q, \partial \overline{\Delta}) = R > M.  
\end{equation}
In particular, $p$ is not on $\partial \Delta$, so there is a 
2-cell whose boundary contains $p$, and hence, by Lemma~\ref{t corridors}, there is a $t$-corridor $\tau$ passing through $p$. 
Applying Lemma~\ref{doubling corridors} to $p$ and $\tau$, we obtain vertices $p_0$ and $p_1$ on the side of $\tau$ not containing $p$, and corridors $\tau_0$ and $\tau_1$ containing $p_0$ and $p_1$ respectively, such that $\tau, \tau_0$ and $\tau_1$ are distinct, and $d(p, p_i) <M$ for $i=0, 1$. 

This process can be iterated as long as Lemma~\ref{doubling corridors} applies --- see Figure~\ref{fig:nested}.  Specifically, suppose that for every $j \le l-1$ and every binary sequence $i_1\cdots i_j$ of length $j>0$, we have constructed, using Lemma~\ref{doubling corridors}, a corridor $\tau_{i_1\cdots i_j}$ which contains a vertex $p_{i_1\cdots i_j}$, such that the corridors are all distinct from each other and $d(p_{i_1\cdots i_j}, p_{i_1\cdots i_{j-1}})<M$. (We understand $p_{i_1\cdots i_{j-1}}$ to be $p$ when $j=1$.)  It follows that 
\begin{equation}\label{pR ineq}
d(p, p_{i_1\cdots i_{l-1}}) \le (l-1)M.
\end{equation}

Now suppose $l<R/M$.  Then $R - (l-1)M >M$, and by the triangle inequality followed by \eqref{pR ineq} and \eqref{pqR ineq}, we have
\begin{equation} 
 d(p_{i_1 \cdots i_{l-1}}, \partial \Delta)  \geq d(p, \partial \Delta) - d(p_{i_1 \cdots i_{l-1}}, p) \geq R - (l-1) M >M. 
 \end{equation} 
Thus, for each binary sequence of length $l-1$, Lemma~\ref{doubling corridors} applies to $\tau_{i_1 \cdots i_{l-1}}$ and $p_{i_1 \cdots i_{l-1}}$, and identifies a pair of distinct new corridors  $\tau_{i_1\cdots i_l}$, for $i_l =0,1$, and a pair of vertices $p_{i_1\cdots i_l}$ on the opposite side of $\tau_{i_1 \cdots i_{l-1}}$, such that $d(p_{i_1\cdots i_l}, p_{i_1\cdots i_{l-1}})<M$.  Given how each corridor separates the diagram into two subdiagrams and how the successive generations nest, all the corridors identified in the course of the iteration must be distinct.  Thus, a total of $2^{l+1}-1$ distinct corridors are found.

 \begin{figure}[ht]
\begin{overpic}
{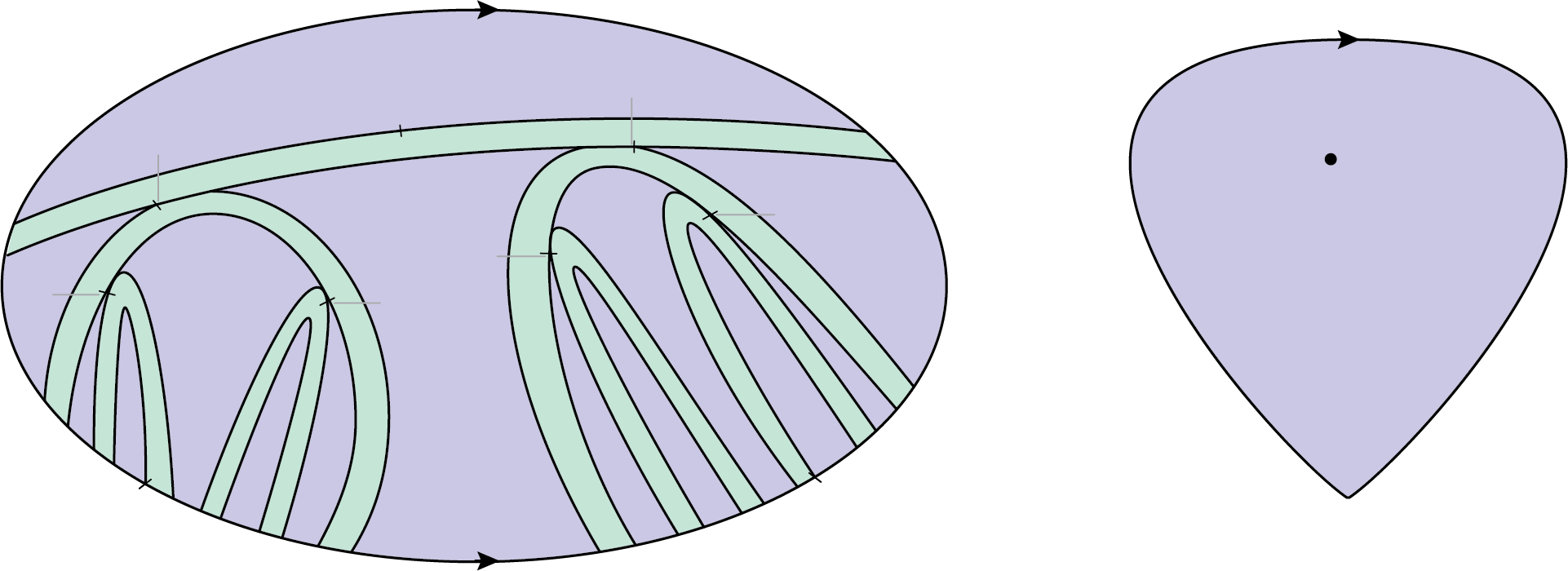}
\put(30,37){\small{$w$}}  
\put(30.3,-2){\small{$v$}}  
\put(7,34){\small{$\Delta$}}  
\put(44,27.2){\small{$\tau$}}      
\put(22.3,13){\small{$\tau_0$}}      
\put(33.5,13){\small{$\tau_1$}}      
\put(10.7,12){\small{$\tau_{00}$}}      
\put(14.2,14.5){\small{$\tau_{01}$}}      
\put(37.4,21){\small{$\tau_{10}$}}      
\put(40.2,19){\small{$\tau_{11}$}}      
\put(24.5,29.5){\small{$p$}}  
\put(9.5,27.5){\small{$p_0$}}      
\put(39.5,30.8){\small{$p_1$}}      
\put(0.5,17){\small{$p_{00}$}}      
\put(25,16.7){\small{$p_{01}$}}      
\put(27.5,19.7){\small{$p_{10}$}}      
\put(50,22.5){\small{$p_{11}$}}     
\put(-1,20.5){\small{$t$}}  
\put(2,8){\small{$t$}}  
\put(5.7,5){\small{$t$}}  
\put(9.5,3){\small{$t$}}  
\put(12.5,1.5){\small{$t$}}  
\put(16.5,0){\small{$t$}}
\put(20.5,-0.5){\small{$t$}}  
\put(39,-1){\small{$t$}}  
\put(44,0.5){\small{$t$}}  
\put(48,2){\small{$t$}}  
\put(52,4){\small{$t$}}  
\put(55.5,6.5){\small{$t$}}  
\put(58,9.5){\small{$t$}}  
\put(57,27){\small{$t$}}  
\put(64,18){\huge{$\to$}}      
\put(85,35){\small{$\overline{w}$}}  
\put(74,34){\small{$\overline{\Delta}$}}  
\put(84.5,27.5){\small{$q$}}  
\end{overpic}
 \caption{Nesting $t$-corridors per our proof of Proposition~\ref{prop lower bound}}
  \label{fig:nested}
\end{figure}

In particular, taking $l = \lfloor R/M \rfloor-1$, we see that $\Delta$ has at least $2^{\lfloor R/M \rfloor}-1$ distinct $t$-corridors, and each of these contributes two distinct edges to $\partial \Delta$.  So $|v| + |w| \ge 2( 2^{\lfloor R/M \rfloor}-1) = 2^{\lfloor R/M \rfloor+1} -2 \ge 2^{R/M}-2$, which proves \eqref{eqn: the aim}.
\end{proof}

\bibliographystyle{alpha}
\bibliography{references}

\end{document}